\documentclass[11pt]{article}

\usepackage[margin=1.15in]{geometry}
\usepackage{amsmath,amssymb,amsthm}
\usepackage{microtype}
\usepackage{hyperref}

\newtheorem{theorem}{Theorem}[section]
\newtheorem{corollary}[theorem]{Corollary}

\newcommand{\E}{\mathbb E}
\newcommand{\F}{\mathbb F}

\title{\bfseries A Counterexample to Talagrand's Operator Cotype Problem}
\author{Xinglong Wu}
\date{}

\begin{document}

\maketitle

\begin{abstract}
Talagrand asked whether the Rademacher cotype of an operator is
universally controlled by the maximum of its Gaussian cotype and its
$(q,1)$-summing norm. We give a negative answer already for $q=2$. The counterexample was discovered by ChatGPT (GPT-5.6).

\end{abstract}

\section{Introduction}

Let $U:X\to Y$ be a bounded operator between Banach spaces. For $q\ge2$,
let $C_q^g(U)$ and $C_q^r(U)$ denote respectively the least constants
such that
\[
\left(\sum_i\|Ux_i\|^q\right)^{1/q}
\le
C_q^g(U)\,
\E\left\|\sum_i g_i x_i\right\|
\]
and
\[
\left(\sum_i\|Ux_i\|^q\right)^{1/q}
\le
C_q^r(U)\,
\E\left\|\sum_i\varepsilon_i x_i\right\|,
\]
for every finite sequence $(x_i)\subset X$. Here $(g_i)$ are independent
standard Gaussians and $(\varepsilon_i)$ are independent symmetric signs.

We also denote by $\|U\|_{q,1}$ the least constant such that
\[
\left(\sum_i\|Ux_i\|^q\right)^{1/q}
\le
\|U\|_{q,1}
\max_{\eta_i=\pm1}
\left\|\sum_i\eta_i x_i\right\|.
\]

Talagrand asked whether there exists a universal constant $L$ such that
\[
C_q^r(U)
\le
L\max\{C_q^g(U),\|U\|_{q,1}\}
\tag{1.1}
\]
for every bounded operator $U$; see
\cite[Research Problem 19.1.2]{Talagrand2021}.
The following gives a negative answer already for $q=2$.

\begin{theorem}\label{thm:main}
There are universal constants $c,C>0$ such that, for every sufficiently
large $n=2^k$, there exist an $n$-dimensional Banach space $X_n$ and an
operator
\[
U_n:X_n\longrightarrow\ell_\infty^n
\]
such that
\[
{C_2^r(U_n)}
\ge c\sqrt{\log(n+1)} {\max\{C_2^g(U_n),\|U_n\|_{2,1}\}}.
\]
\end{theorem}
For a Banach space \(X\), let \({\mathrm{Kw}}(X)\) denote the least
constant \(K\) such that every finite sequence \(x_1,\ldots,x_m\in X\)
admits a decomposition
\[
x_i=x_i'+x_i''
\]
satisfying
\[
\E\left\|\sum_i g_i x_i'\right\|
\le
K\E\left\|\sum_i\varepsilon_i x_i\right\|
\]
and
\[
\max_{\eta_i=\pm1}
\left\|\sum_i\eta_i x_i''\right\|
\le
K\E\left\|\sum_i\varepsilon_i x_i\right\|.
\]
Kwapień's decomposition problem asks whether this constant can be chosen independently of the Banach space; equivalently, whether there exists a universal constant \(K<\infty\) such that
\[\mathrm{Kw}(X)\le K\] holds for all Banach spaces $X$.

The problem was recorded in its infinite-series form by Bednorz and Latała as Problem 1.3 in \cite{BednorzLatala2014}, where it was attributed to S. Kwapień, and was later stated in the above finite-sequence form by Talagrand as Research Problem 19.1.3 in \cite{Talagrand2021}.
As observed by Talagrand
\cite[Exercise 19.1.4]{Talagrand2021}, such a decomposition implies
\[
C_2^r(U)
\le
2\mathrm{Kw}(X)
\max\{C_2^g(U),\|U\|_{2,1}\}
\]
for every bounded operator \(U:X\to Y\).
Theorem \ref{thm:main} immediately yields a negative answer for Kwapień's  decomposition problem.
\begin{corollary}\label{cor:kwapien}
There exists a universal constant \(c>0\) such that, for every
sufficiently large \(n=2^k\), the spaces \(X_n\) in
Theorem~\ref{thm:main} satisfy
\[
\mathrm{Kw}(X_n)
\ge
c\sqrt{\log(n+1)}.
\]
 
\end{corollary}
\subsection*{Statement on A.I. Use}
The counterexample was discovered during discussions with ChatGPT (GPT-5.6). The author has independently verified all arguments presented in the paper.
\section{Proof of theorem 1.1}

Fix $n=2^k$ and let $H=H_n$ be the normalized Walsh--Hadamard matrix,
indexed by $\F_2^k$:
\[
H_{ab}=n^{-1/2}(-1)^{a\cdot b}.
\]
Thus
\[
H^TH=I_n,
\qquad
|H_{ab}|=n^{-1/2}.
\]
This is the standard Sylvester–Walsh construction; see, forexample, \cite{Sylvester1867} for its classical origin.
Set
\[
s_n=\sqrt{2\log(2n)}
\]
and equip $\mathbb R^n$ with the norm
\[
\|x\|_{X_n}
=
\max\left\{
\|x\|_\infty,\,
\frac{\|Hx\|_\infty}{s_n}
\right\}.
\]
Let
\[
U_n:X_n\longrightarrow\ell_\infty^n,
\qquad
U_nx=x.
\]

We estimate the three quantities appearing in
Theorem~\ref{thm:main}.

First, if $\varepsilon=(\varepsilon_1,\ldots,\varepsilon_n)$ is a sign
vector, then every coordinate of $H\varepsilon$ is $1$-subgaussian.
Hence
\[
\E\|H\varepsilon\|_\infty
\le
\sqrt{2\log(2n)}
=s_n.
\]
Since $\|\varepsilon\|_\infty=1$,
\[
\E\|\varepsilon\|_{X_n}
\le
1+\frac{\E\|H\varepsilon\|_\infty}{s_n}
\le2.
\]
Applying the definition of Rademacher cotype to the standard basis gives
\[
C_2^r(U_n)
\ge
\frac{\sqrt n}{2}.
\tag{2.1}
\]

Second, the classical estimate
\[
C_2^g(I_{\ell_\infty^n})
\le
C\sqrt{\frac{n}{\log(n+1)}}
\tag{2.2}
\]
holds; see, for example, \cite{GeissJunge1995}. Since
$\|x\|_\infty\le\|x\|_{X_n}$, it follows immediately that
\[
C_2^g(U_n)
\le
C\sqrt{\frac{n}{\log(n+1)}}.
\tag{2.3}
\]

It remains to estimate the $(2,1)$-summing norm. Let
$x_1,\ldots,x_m\in X_n$, write
\[
A=[x_1\ \cdots\ x_m]=(a_{ji}),
\qquad
HA=(c_{ji}),
\]
and set
\[
\alpha=\max_j\sum_i|a_{ji}|,
\qquad
\beta=\max_j\sum_i|c_{ji}|.
\]
Then
\[
D:=
\max_{\eta_i=\pm1}
\left\|\sum_i\eta_i x_i\right\|_{X_n}
=
\max\left\{\alpha,\frac{\beta}{s_n}\right\}.
\tag{2.4}
\]

Put $u_i=\|x_i\|_\infty$. Clearly $u_i\le\alpha$. Since
$A=H^T(HA)$ and $|H_{jr}|=n^{-1/2}$,
\[
u_i
\le
\frac1{\sqrt n}\sum_{r=1}^n|c_{ri}|.
\]
Therefore
\[
\sum_i u_i^2
\le
\frac{\alpha}{\sqrt n}
\sum_{r,i}|c_{ri}|
\le
\sqrt n\,\alpha\beta.
\]
Using \emph{(2.4)},
\[
\left(\sum_i\|U_nx_i\|_\infty^2\right)^{1/2}
\le
n^{1/4}(\alpha\beta)^{1/2}
\le
n^{1/4}s_n^{1/2}D.
\]
Hence
\[
\|U_n\|_{2,1}
\le
n^{1/4}(2\log(2n))^{1/4}.
\tag{2.5}
\]

Combining \emph{(2.1)}, \emph{(2.3)}, and \emph{(2.5)}, and observing that
\[
n^{1/4}(\log n)^{1/4}
=
o\left(\sqrt{\frac n{\log n}}\right),
\]
we obtain, for all sufficiently large $n$,
\[
\max\{C_2^g(U_n),\|U_n\|_{2,1}\}
\le
C\sqrt{\frac n{\log(n+1)}}.
\]
Consequently,
\[
\frac{C_2^r(U_n)}
{\max\{C_2^g(U_n),\|U_n\|_{2,1}\}}
\ge
c\sqrt{\log(n+1)},
\]
which proves Theorem~\ref{thm:main}.

\end{document}